\documentstyle[12pt]{article}

\makeatletter
\def\@begintheorem#1#2{\it \trivlist \item[\hskip
\labelsep{\bf #2\ #1.}]}
\def\@opargbegintheorem#1#2#3{\it \trivlist
      \item[\hskip \labelsep{\bf #1\ #2\ (#3).}]}
\makeatother

\newtheorem{them}{Theorem}[section]

\newtheorem{lemm}[them]{Lemma}
\newtheorem{cor}[them]{Corollary}

\newtheorem{thm}[them]{Theorem}

\newcommand{\jtbnumpar}[1]{\refstepcounter{them}
\trivlist
\item[\hskip \labelsep{\bf \thethm \ #1.}]}

\newcommand{\jtbdef}{\jtbnumpar{Definition}}
\def\jtbnot{\jtbnumpar{Notation}}

\def\PROOF #1.{\par\noindent{\it Proof#1}.\ \ignorespaces}
\def\proof #1.{\par\noindent{\it Proof#1}.\ \ignorespaces}

\def\subm{\leq}

\def\Bscr{{\cal B}}

\def\Bscr{{\cal B}}

\let\?=\joinrel

\let\leftv=^
\let\rightv=^

\def\k/{\kern.2em}    

\def\rng{\mathop{\rm rng}}
\def\im{\mathop{\rm rng}}

\def\cl{\mathop{\rm cl}}

\def\max{\mathop{\rm max}}

\let\union=\cup             %

\edef\bigcup{\mathop{\textstyle\mathchar\the\bigcup}}
\let\bigunion =\bigcup

\let\inter=\cap             %

\edef\bigcap{\mathop{\textstyle\mathchar\the\bigcap}}

\edef\bigwedge{\mathop{\textstyle\mathchar\the\bigwedge}}

\edef\bigvee{\mathop{\textstyle\mathchar\the\bigvee}}

\edef\sum{\mathop{\textstyle\mathchar\the\sum}}
\def\math&{\ \& \ }

\def\force {\mathrel^\joinrel\rightarrow}
\def\force {\mathrel{\scriptstyle\mathrel^\joinrel\rightarrow}}
\def\forceq {\mathrel{\mathop{\force}\limits_{\textstyle\texsim}}}
\def\forceq{\mathrel^\joinrel
 \mathrel{\mathop{\rightarrow}\limits_{\smash{\textstyle\texsim}}}}
\def\forceq{\mathrel{\scriptstyle\mathrel^\joinrel
 \mathrel{\smash{\mathop{\rightarrow}\limits_{\smash{\raise
 2pt\hbox{$\scriptstyle\texsim$}}}}}}}

\let\exclaim=!                 %
\let\iso\approx             %
\let\texsim=\sim         %
\let\conj\sim             %
\def\conjp #1 {\conj_{#1}}     %
\let\sim\simeq             %
\let\neg=\lnot             %

\def\0bar{\bar 0}         %
\def\1bar{\bar 1}         %
\def\abar{\overline a}

\let\sat=\models                %

\def\bbar{\overline b}
\def\cbar{\overline c}

\def\hbar{\overline h}

\def\wbar{\overline w}
\def\xbar{\overline x}
\def\ybar{\overline y}
\def\zbar{\overline z}

\def\Proof{Proof}

\def\bigunion{\union}

\def\Dscr{{\cal D}}

\def\text#1{\ifmmode\leavevmode\hbox{#1}\else
   \typeout{Warning: \string\text \space used outside math mode!}
   \begingroup\hbox{#1}\endgroup\fi}

\def\bar#1{\overline{#1}}
\def\implies{\rightarrow}

\def\subi{\leq_{i}}
\def\subm{\leq_{s}}

\def\bK{{\bf K}}
\def\Var{{\rm Var}}
\def\Cov{{\rm Cov}}

\title{Randomness and Semigenericity\footnote{Mathematics Review  
Numbers:
03C10, 05C80, keywords: random graphs, $0-1$-laws, stability}}
\author{John T. Baldwin
\thanks{Partially supported by NSF grant
9308768 and a visit to Simon Fraser University}
\\ Department of Mathematics, Statistics and  Computer
Science\\University of Illinois at Chicago
\and
Saharon Shelah\thanks{This is paper 528.  Both authors thank
Rutgers University and the Binational Science Foundation for partial
support of this research.}\\ Department of Mathematics\\
Hebrew University of Jerusalem\\Rutgers University}

\begin{document}
\maketitle
\begin{abstract}
Let $L$ contain only the equality symbol and let
$L^+$ be an arbitrary finite symmetric relational language containing  
$L$.
Suppose
probabilities are defined on finite $L^+$ structures with \lq edge  
probability"' $n^{-\alpha}$. By $T^{\alpha}$, the almost sure theory  
of random $L^+$-structures
we mean the collection of $L^+$-sentences which have limit
probability
1.     $T_{\alpha}$ denotes the theory of the generic
structures
for $\bK_{\alpha}$, (the collection of finite graphs $G$ with
$\delta_{\alpha}(G) =|G| - \alpha \cdot |\text{ edges of G }|$  
hereditarily nonnegative.)

\begin{thm}  

$T^{\alpha}$, the almost sure theory of random
$L^+$-structures is the same as the theory $T_{\alpha}$
of the $\bK_{\alpha}$-generic model.
This theory is
complete, stable, and nearly model complete.  Moreover, it has the
finite model property and has only infinite models
so is not finitely
axiomatizable.
\end{thm}
\end{abstract}

This paper unites two apparently disparate lines of
research.  In \cite{ShelahSpencer}, Shelah and Spencer proved a
$0-1$-law
for first order sentences about random graphs with edge probability
$n^{-\alpha}$ where  $\alpha$ is an irrational number between $0$ and
$1$.  Answering a question raised by Lynch \cite{Lynch1},
we extend this result from graphs to hypergraphs (i.e. to arbitrary
finite symmetric relational languages).
Let $T^{\alpha}$
denote the set of sentences with limit probability $1$.  The
Spencer-Shelah proof
proceeded by a process of quantifier elimination which implicitly
showed the
theories $T^{\alpha}$ were {\em nearly model complete} (see below)
and complete.

Hrushovski in \cite{Hrustableplane} refuted a conjecture of Lachlan
by constructing an $\aleph_0$-categorical strictly stable
pseudoplane.
Baldwin and
Shi \cite{BaldwinShiJapan} considered a variant on his methods to
construct  strictly stable (but not
$\aleph_0$-categorical) theories $T_{\alpha}$ indexed by irrational
$\alpha$.
In this paper we
show that for each irrational
$\alpha$, $T^{\alpha} = T_{\alpha}$
and thus deduce that $T_{\alpha}$ is not
finitely axiomatizable and that $T^{\alpha}$ is stable.

Each $T_{\alpha}$ is the theory of a `generic' model
$M_{\alpha}$ of an amalgamation class $\bK_{\alpha}$ of finite
structures.  Although the Hrushovski examples are easily seen to be
nearly model complete this is less clear for the $T_{\alpha}$ since
they are not
$\aleph_0$-categorical.
We show that each
$T_{\alpha}$ is nearly model complete.

In the first, purely model theoretic,
section of the paper we describe our basic framework and
prove a sufficient condition for certain theories, including the
$T_{\alpha}$,
to be nearly model
complete.  These conditions depend upon a generalization of the
notion of  {\em genericity} of a structure: {\em semigenericity},  
which is
introduced
in this paper.
In the second section we consider the addition of  random relations
and deduce the main results for this case:
The almost-sure theory and the theory of the generic model are equal,
complete, stable, nearly model complete, and not finitely
axiomatizable.
\relax From the model theoretic standpoint the extension from graphs to an
arbitrary finite relational language is not a big step; it was 

spelled out in \cite{Wagnerdim}.   The distance is larger from the  
probability standpoint and the problem of making such an extension  
had
been raised by Lynch \cite{Lynch1}.

The first author greatly benefited from discussions on this paper
with
M. Albert, G. Cherlin, M. Itai,
A.H. Lachlan, C. Laskowski, D. Kueker and D.
Marker.  We want to thank Shmuel Lifsches for a careful reading of
Section~\ref{nmc}.

\section{Near model completeness}
\label{nmc}

After Hrushovski's construction of counterexamples to the conjectures
of Lachlan and
Zil'ber a number of authors explored generalizations of the
variation he had introduced on the Fraiss\'e-Jonsson construction.
Hrushovski had noted that in his situation, where the generic model
was $\omega$-saturated, the theory of the generic admitted the level
of quantifier elimination which we christen "nearly model complete"  
in 

this paper.
We
reprise one general setting for this study here and in the next
section connect it with certain random models.  Baldwin and Shi
\cite{BaldwinShiJapan} studied a situation where the
homogeneous-universal model, renamed {\em generic} by Kueker and  
Laskowski
\cite{KuekerLas}, is not $\omega$-saturated.  Kueker and Laskowski
investigated the conditions in which the theory constructed from a  
generic admitted
various levels of quantifier elimination.  After the first author
noticed the connection between \cite{Hrustableplane} and
\cite{ShelahSpencer}, we began to consider the quantifier complexity
of the theory $T^{\alpha}$. There is no explicit elimination of  
quantifiers
result in \cite{ShelahSpencer} but a lemma similar to our  
Lemma~\ref{heart}
which is the crucial technical step.
 The second author had already begun
notes
generalizing \cite{ShelahSpencer}; the $0-1$-law in
Section~\ref{addran} contains a more concrete version of his
approach.
Shelah has continued this approach to the probability aspect in more
generality in \cite{Shelah550}.  A close look at the quantifier
elimination results in \cite{Hrustableplane,ShelahSpencer}, suggests  
the following
definition.

\jtbdef A theory $T$ is said to be {\em nearly model complete} if
every formula is equivalent in $T$ to a Boolean combination of
$\Sigma_1$-formulas.

Thus, $T$ is nearly model complete if the type of any finite sequence
is determined by exactly the family of $\Sigma_1$-formulas it
satisfies.  Near model completeness lies strictly in strength between
model completeness and $1$-model completeness (every formula is
equivalent to
a $\Sigma_2$-formula).

\jtbnot Fix a finite relational language $L$.  For any class $\bK$ of
structures, $S(\bK)$ denotes the class of all substructures of
members
of $\bK$.  Let $\bK_0$ be a collection of finite $L$-structures and
{\bK} be a class of models whose finite substructures are in
$S({\bK_0})$.  We always assume that the empty structure is in
$\bK_0$.
We will consider several different choices for $\bK_0$
in this paper.  In the following, $A$, $B$, $C$ vary over
$\bK_0$; $M$, $N$ over $\bK$.  If $A$, $B$ are subsets of $N$, we
write $AB$ for the $L$-structure contained in $N$ with universe $A
\union B$.

If $B \inter C = A$ we write
$B\otimes_A C$
for the structure with
universe $B\union C$ and no relations other than those on $B$ or $C$.
If $A,B,C$ are substructures of $N$ such that the structure imposed
by
$N$ on $BC$ is isomorphic to
$B\otimes_A C$ we say $B$ and $C$ are {\em freely joined over $A$ in
$N$}.
In general we do not assume $\bK_0$ is closed under $\otimes$ but
this
assertion will turn out to be an important property of some classes
we
consider.
We write $X\subseteq_{\omega} Y$ to indicate $X$ is a finite subset
of
$Y$.

\vskip
.2in We will first discuss a class of finite structures equipped with
a dimension function satisfying certain natural properties.  Then we
define from this dimension function a notion of strong submodel.  The
main quantifier elimination result is proved in terms of the strong
submodel concept.  But, the connection with random models is
obtained by exploiting an appropriate dimension function.  The fact
that
 this dimension function (in Example~\ref{kalpha}) is the same as
that
employed by Hrushovski to construct a strictly
stable $\aleph_0$-categorical
pseudoplane is the key to the argument for the stability of the
almost
sure theory of random graphs with edge probability $n^{-\alpha}$ when
$\alpha$ is irrational.

\jtbdef
\label{setup}
Let $\delta$  be an
arbitrary function assigning a real number to each isomorphism type
of
finite $L$-structure with $\delta(\emptyset) = 0$.  $\delta(A/B)$
equals
by definition $\delta(AB)-\delta(B)$.  This yields immediately:
$$\delta(AB/C) = \delta(A/BC) + \delta(B/C).$$

Note that the structure with universe $AB$ (and thus $\delta(A/B)$)
is
not determined by the separate structures on $A$ and $B$ but by some
embedding of both into an element of $\bK_0$.

\jtbnot
\label{fixnotation}
  We deal only with
structures on which the relations of $L$ are symmetric (i.e
$R(\abar)$
holds just if it holds for any permutation of $\abar$) and  
irreflexive
(i.e. hold only for
sequences of distinct elements).  Thus the relations are on sets
rather than sequences.

 We require the following conditions on $\delta$.

\jtbnumpar{Axiom}   $\bK_0$ and $\delta$ satisfy for $A,B,C \ldots
\in S(\bK_0)$ and $N,M, \in \bK$:

\label{yax}
\begin{enumerate}
\item  $\delta: S({\bf K}_0) \mapsto \Re^+$ (the nonnegative reals)
and $\delta(\emptyset) = 0$.
\item If $A$, $B$, and $C$ are disjoint subsets of $N$ then
  $\delta(A/B) \geq \delta(A/BC)$.
\item For every $n\in \omega$ there is an $\epsilon_n>0$ such that if
  $|C|<n$ and $A, C$ are disjoint subsets of $M$ with $\delta(CA/A) <
  0$ then $\delta(CA/A) \leq -\epsilon_n$.
\item There is a real number $\epsilon>0$ such that if $A,B,B'$ are
  disjoint subsets of a model $N$ and $\delta(A/B) - \delta(A/BB') <
  \epsilon$ then $R(A,B,B') = \emptyset$ and $\delta(A/B) =
  \delta(A/BB')$.
\item If $f$ is a $1-1$ homomorphism from $A$ to $B$ then for every
$X \subseteq Y \subseteq A$, $\delta(Y/X) \geq \delta(f(Y)/f(X))$.
\end{enumerate}

Axioms iii) and iv) play no explicit role in the argument presented
here.  But they are important in establishing the stability of
$T_{\alpha}$ in \cite{BaldwinShiJapan} so are used in the proof of
Theorem~\ref{talphaexists}.  Note that Axiom~\ref{yax}~iv) is
stronger
than the assertion that if $f$ is a $1-1$-homomorphism,
$\delta(X) \geq \delta(f(X))$.

Axiom~\ref{yax} i) requires that the range of $\delta$ be the
nonnegative
reals.  This allows us to obtain an
important monotonicity property by modifying $\delta$ to $d: {\bK}
\times
S({\bf K}_0) \mapsto \Re^+$ by defining for each $N \in {\bK}$,
$$d(N,A) = \inf
\{\delta(B): A \subseteq B \subseteq_{\omega} N\}.$$
We usually write
$d(N,A)$ as $d_N(A)$.  We will omit the subscript $N$ if it is clear
from context.  This operator serves only as a notational convenience
within this paper
but plays an essential role in establishing the stability of
$T_{\alpha}$ in
\cite{BaldwinShiJapan}.  The nonnegativity requirement on $\delta
|\bK_0 $
not only justifies the definition of $d_N(A)$  but is necessary for
the important Lemma~\ref{h1}.

The classes $({\bf K_{\alpha}, \delta_{\alpha}})$, which are
defined as follows, are important examples of this situation.

\jtbnumpar{Example}
\label{kalpha}
Let the relation symbols of $L$ be $\langle R_i: i < p\rangle$.  Let
$w_i(A)$ be the
 realizations of $R_i$ in $A$.  Fix a sequence
$\bar{\alpha}$ with
 $0< \alpha_i \leq 1$ for $i<p$.  Then for each $A$, let $e(A) = \sum  
w_i
\alpha_i$.
Let $\bK_{\alpha}$ denote the class of all finite $L$-structures $A$
such that
for all substructures $A'$ of $A$, $$\delta_\alpha(A') = |A'| -
e(A') \geq 0.$$
See \cite{BaldwinShiJapan} for
the straightforward verification of the axioms in this example.

\medskip From the dimension function we define certain
special notions of submodel which make it easier to formulate
our argument.

\jtbdef For finite $A$, $B$ contained in $N$, define the relative
dimension of $A$ over $B$, $d_N(A/B)$ as $d_N(A/B) = d_N(AB) -
d_N(B)$.  If $B
\subseteq A \subseteq_{\omega} N $ this simplifies to $d_N(A/B)
=d_N(A) -
d_N(B)$

\jtbdef
\label{subdef}
For $M \subseteq N \in S({\bK})$, define $M\subm N$ if for each
finite
$X \subseteq M$, $d_M(X) = d_N(X)$.  We say $M$ is a {\em strong}
submodel of $N$.
We say $f:M \mapsto N$ is a {\em strong embedding} if $fM \subm N$.
We write $M <_s N$ if $M \subm N$ but $M \neq N$.

\smallskip

We introduce a second kind of distinguished substructure by defining
$\subi$ from $\subm$ as follows.  Note that the definition yields
that
$A \subi A$.
\jtbdef
\label{defins1}
For $A,B \in S(\bK_0)$, $A \subi B$ if $A \subseteq B$ but there is
no
$A'$ with $A \subseteq A' <_s B$.  If $A \subi B$, we say $B$ is an
{\em intrinsic} extension of $A$.

In terms of the dimension function $A \subi B$ means $A=B$ or
$\delta(B/A) < 0$ and $\delta(B/A) < \delta(B'/A)$ for any
intermediate $B'$.

\begin{lemm}
\label{yaxtrue}  Consider the situation described in
Definition~\ref{setup}.  If $\delta$ is a dimension function
satisfying the properties of Axiom~\ref{yax} and $\subm$ is defined
as
in Definition~\ref{subdef} then $( \bK, \subm)$ satisfies the
following
conditions for $M,N,N'\in S(\bK)$.

{ \bf A1}.  $M\subm M$.

{ \bf A2}. If $M\subm N$ then $M\subseteq N$.

{ \bf A3}. If $M\subm N\subm N$ then $ M\subm N'$.

{ \bf A4}.
If $M\subm N$, $N'\subseteq N$
then $M\bigcap N' \subm N'$.

{\bf A5}.  $\subi$ is preserved under $1-1$ homomorphism.

{\bf { A6}.} For all $M\in S(\bK)$, $\emptyset \subm M$.
\end{lemm}

\jtbnumpar{Remark} Note that { \bf A4} implies that if $M,N,N'\in
S(\bK)$ and $M \subm N$, $N' \subseteq N$ and $M\subseteq N'$ then $M
\subm N'$.

The quantifier elimination results of this section could be obtained  
by
taking as primitive a class
$\bK_0$ equipped with a notion of strong submodel, and regarding the
results of
Lemma~\ref{yaxtrue} and Lemma~\ref{h1} below as axioms. Naturally, we  
would
then require that
$\bK_0$ and $\subm$ be closed under isomorphism.  The dimension
function is needed for the calculations in Section 2.

\jtbnumpar{Remark} In earlier formulations, the relation $\subm$ was
defined just on $\bK_0$ rather than on $S(\bK_0)$.  This leads to
difficulties in phrasing Axiom {\bf A4}.  Our current formulation
extends the ideas of \cite{BaldwinShiJapan} to encompass the
Baudisch construction of a new $\aleph_1$-categorical group
\cite{Bgroup}.  For our purposes in this paper, we could have
identified $\bK_0$ with $S(\bK_0)$ and we make that restriction in
Section 2.

\jtbnumpar{Remark}
\label{fixingA6
}
Axiom {\bf A6} holds in the examples at hand because the range of
$\delta$ is nonnegative as specified in Axiom~\ref{yax}.
In Section 2 we will begin with a $\delta$ mapping all finite
$L$-
structures into the reals.  The requirement that $\delta$ is
nonnegative  requires revising the choice of
$\bK_0$ (and thus $\bK$) to guarantee that if $M \in \bK$, then for
every $A
\subseteq_{\omega} M$, $\delta(A) \geq 0$.  We show it is harmless
to make this assumption in Lemma~\ref{justa6}.

\smallskip

\begin{lemm}{}
\label{dichot1}
\begin{enumerate}
\item For any $A \subseteq C\in \bK_0$, we can choose $B$ with $A
\subi B \subm
  C$.
\item
$\subi$ is transitive.
\end{enumerate}
\end{lemm}

\Proof.  For i), let $B$ have minimal cardinality among the subsets
$X$ of $C$ that contain $A$ with $X \subm C$.  Use {\bf A4} for ii).

\jtbdef For any $L$-structure $M$, let $A,B$ be finite substructures
of $M$ with $A \subseteq B$.
Then
\begin{enumerate}
\item By a {\em copy} of $B$ over $A$ in $M$ we mean the image of an
extension $\hat f$
to $B$
of an embedding $f$ from $A$ into $M$.
\item $\chi_M(B/A)$ is the number of distinct copies of $B$ over $A$
  in $M$.
\item
$\chi^*_M(B/A)$
is the supremum of the cardinalities of maximal
  families of disjoint (over $A$) copies of $B$ over $A$ in $M$.
\end{enumerate}

\begin{lemm}  If $\Bscr$ and $\Bscr'$ are maximal
  families of disjoint over $A$ copies of $B$ over $A$ then $|\Bscr|
  \leq |B-A| |\Bscr'|$.
\end{lemm}

Proof.  Define a map from $\Bscr$ to $\Bscr'$ by mapping each element
of $\Bscr$ to an element of $\Bscr'$ that it intersects off $A$.
This
map is at most $|B-A|$-to-one since the members of $\Bscr$ are
disjoint over $A$.

In particular, this shows that the supremum in the definition of
$\chi^*_M$
is achieved.
As one varies over the entire family of examples of structures
constructed in this manner (e.g. in \cite{Hrustableplane},
\cite{BaldwinShiJapan}, etc.)
the dimension function produces an important trichotomy concerning
pairs $A \subset B$.  Consider an infinite $(\bK_0,\subm)$-generic
(Definition~\ref{gendef}) model $M$.
$\chi^*_M(B/A)$ will be bounded if $A\subi B$, infinite if
$\delta(B/A) >
0$, and will vary with the choice of $(\bK_0,\subm)$ 

if $\delta(B/A) =0$.  The key to the
$0-1$ law in Section~\ref{addran} is that when ${\alpha}$ denotes
a sequence, which is linearly independent with $1$,
the third case cannot occur.  The uniform bound on $\chi^*_M(B/A)$
follows from
our restricting $\bK_0$ so $\delta$ is nonnegative.  In
\cite{Shelah467},
Shelah
considers  a different probability measure which does not permit the
nonnegativity restriction; in that situation $\chi^*_M(B/A)$ is a  
slow
growing function.
In our situation we have the following.

\begin{lemm}
\label{h1}There is a binary function $t:\omega \times \omega \mapsto
\omega$
which is monotone increasing in both arguments such that if $A \subi
B$ then for any $M \in \bK$ with $A \subseteq
  M$, $\chi_M(B/A) \leq t(|A|,|B|)$.
\end{lemm}

\jtbnumpar{Remark}  This follows easily from Lemma 3.19 of
\cite{BaldwinShiJapan}.  One must note
that $A \subi B$ if and only if there is a sequence $A = A_0, A_1,
\ldots A_n = B$ such that
$(A_i,A_{i+1})$ is a minimal pair in the sense of
\cite{BaldwinShiJapan}.

\jtbdef
 For any $M \in \bK$, any $m\in \omega$, and any $A \subseteq M$,
  $$\textstyle{\cl^m_M(A)} = \bigunion\{B:A \subi B \subseteq M \ \&
 \  |B-A| < m\}.$$

\medskip The following are immediate from Lemma~\ref{h1} and the
definitions.

\begin{lemm}
\label{bound}
There is a function $f$ mapping $\omega \times \omega$ into $\omega$
such that for any $A\subseteq_{\omega} M\in \bK$ and $m$,
$\cl^m_M(A)$
is finite and its cardinality is uniformly bounded by $f(|A|,m)$.
\end{lemm}

\begin{lemm}
\label{compress}
For any $M,A,m,n$ there exists a $p$ depending on $|A|,m,n$ but not
on
the embedding of $A$ into $M$ with $\cl^m_M(\cl^n_M(A)) \subseteq
\cl^p_M(A)$.
\end{lemm}

\proof. Let $p = m + f(|A|,n)$ and check.


The next result is immediate noting that $A \subi X$ does not
depend on any ambient model containing $X$.
\begin{lemm}
\label{1.9a}
If $A \subseteq B \subseteq C$ and $\cl^m_C(A) \subseteq B$, then
$\cl^m_C(A) = \cl^m_B(A)$.
\end{lemm}

\jtbdef
\label{gendef}
 The countable model $M\in \bK$ is {\em $({\bK_0},
  \subm)$-generic} if
\begin{enumerate}
\item If $A\subm M, A\subm B\in \bK_{0}$, then there exists $B'\subm
M
  $ with $ B\cong_{A} B'$; and
\item $M$ is the union of $\langle A_i: i < \omega\rangle$ where each
  $A_i \in \bK_0$ and $A_i \subm A_{i+1}$.
\end{enumerate}

\jtbdef A class $(\bK,\subm)$ has the {\em amalgamation property} if
for any three structures
$A,B,C \in \bK$ with strong embeddings $f,g$ from $A$ into $B$, $C$
there exists $D \in \bK$ and strong embeddings $f': B\mapsto D$,
$g':C\mapsto D$ with $f'f = g'g$.

Following the Fraiss\'e-Jonsson construction, it is easy to show the
following result.
\jtbnumpar{Fact}
\label{genexist}
If $(\bK_0,\subm)$ satisfies {\bf A0} through {\bf A6} and has the
amalgamation property then there is a unique
  countable $(\bK_0,\subm)$-generic model.

\medskip

We need a more local notion.  This is the key new idea of this paper;
it arose from the notion of a full model in \cite{BaldwinShiJapan}
and
from considering the role of $\cl^m_M(A)$ in \cite{ShelahSpencer}.

\jtbdef The countable model $M$ is {\em $(\bK_0,\subm)$-semigeneric},
or
just semigeneric, if
\begin{enumerate}
\item $M \in \bK$
\item If $ A\subm B\in \bK_{0}$ and $g:A \mapsto M$, then for each
  finite $m$ there exists an embedding $\hat g$ of $B$ into $M$ which
  extends $g$ such that
\begin{enumerate}
\item $\cl^{m}_M({\hat g} B) = {\hat g}B \union \cl^{m}(gA)$
\item $M |\cl^m_M(gA){\hat g} B$ is the free join over $gA$ of
$\cl^m_M(gA)$ and ${\hat g}B$.
 \end{enumerate}
 \end{enumerate}

 In our applications any generic model is semigeneric
 (Lemma~\ref{genissemigen}), so Fact~\ref{genexist} provides us with
a
 semigeneric model.  But while generic models are unique there are
 many semigeneric models in the situations that we deal with here.

 We describe below an infinite set of first order formulas
$\phi^m_{A,B,C}$
which  allow us to
 axiomatize the class of semigeneric models by the following lemma,
which is immediate once we have made the definitions. Note that
these are $\Pi_3$-formulas as there is a universal quantifier
hidden
in the last clause.

\begin{lemm}

\label{axsemigeneric}
The structure $N\in \bK$ is
  semigeneric, if and only if for each $A \subm B$,
each $m < \omega$,
and $C\in \Dscr^m_A$,
$N \sat \phi^m_{A,B,C}$
\end{lemm}
In establishing the following notation we are suppressing a fixed
correspondence between
enumerations of the structures $A$, $B$, $C$, $D$ and the variables
$\xbar$, $\zbar$, $\ybar$, $\wbar$.  This correspondence is chosen
to preserve natural inclusions among the structures and the
variables.  Intuitively, $\Dscr^m_A$ is the set of possible
isomorphism types for
$\cl^m_M(A)$.
The structure
$M$ satisfies $\phi^m_{A,B,C}$ just if the
definition of semigenericity holds for the  finite structures $A  
\subseteq
C$
and $A \subm B$ when $C \iso \cl^m_M(A)$.   

 \jtbnot
\begin{enumerate}
\item   Write $A \subi^m D$ if for each $d \in D-A$, there is a $B$
with $Ad \subseteq B$, $A \subi B$, and $|B-A| < m$.
\item
Let $A \in \bK_0$.
$$\Dscr^m_A =
\{D\in \bK_0: A \subi^m D  \}.$$
Note that by Lemma~\ref{bound} if $D \in  \Dscr^m_A$, $|D|<f(|A|,m)
$.
\item For $C \in \Dscr^m_A$ let $\Dscr^m_{A,C }$ be the set of
 $D \in \Dscr^m_A$ which cannot be embedded in $C$.

\item For any finite $A$, $\delta_A(\xbar)$ denotes the atomic
diagram of
$A$.

\item For $C \in \Dscr^m_A$, $\theta^m_{A,C}(\xbar,\ybar)$ is the
formula
$$\delta_A(\xbar) \wedge \delta_C(\xbar,\ybar) \wedge (\forall
\wbar_D) \bigwedge_{D\in
  \Dscr^m_{A,C}}\neg \delta_D(\xbar,\wbar_D).$$  Then for
$\abar\cbar$ an
enumeration of $C$, $A \subseteq C \subseteq N$, $N\sat
\theta^m_{A,C}(\abar,\cbar)$ if and only if $C = \cl^m_N(A)$.

\item For $A \subm B$ and $C\in \Dscr^{m}_A$, let
$\phi^m_{A,B,C}$
be the sentence

$$(\forall \xbar) (\forall \ybar) (\exists \zbar)[ \delta_A(\xbar)
\wedge\theta^m_{A,C}(\xbar,\ybar)\implies (
\delta_{C\otimes_A   B}
(\xbar,\ybar,\zbar) \wedge \theta^m_{B,
{C\otimes_A   B}
}(\xbar,\ybar,\zbar))].$$
\end{enumerate}

\begin{thm}
\label{qethm}
  If $(\bK_0,\subm)$ satisfies {\bf A1-A6} of Lemma~\ref{yaxtrue} and
Lemma~\ref{h1},
  then for every formula $\phi(\xbar)$ there is a Boolean combination
  of existential formulas $\psi_{\phi}(\xbar)$ such that if $M$ is
  $(\bK_0,\subm)$-semigeneric then $\psi_{\phi}(\xbar)$ is equivalent
  to $\phi(\xbar)$ on $M$.
\end{thm}

\Proof.  We first show:

\begin{lemm}
\label{heart}
For any formula $\phi(x_1 \ldots x_r)$ there is an integer $\ell =
\ell_{\phi}$, such that for any pair of semigenerics $M, M'\in \bK$
and
any $r$-tuples $\abar\in M$ and $\abar'\in M'$ if
$\cl^{\ell_{\phi}}_M(\abar) \iso \cl^{\ell_{\phi}}_{M'}(\abar')$ by
an isomorphism taking
$\abar$ to $\abar'$,  then
$M\sat \phi(\abar)$ if and only if $M'\sat \phi(\abar')$.
\end{lemm}

\Proof.  The proof is by induction on formula complexity; the cases
involving
Boolean
connectives are easy.  So suppose $\phi(\xbar)$ is of the form
$(\exists y)
\psi(\xbar,y)$.  Suppose $M \sat \phi(\abar)$, so there is a $b$ such
that $M \sat \psi(\abar,b)$.

Choose $p_1$ large enough so that for any $N\in \bK$, any $r$-tuple
$\cbar$ from $N$ and any $d\in N$, $|\cl^{\ell_\psi}_N(\cbar,d)|<
p_1$.  Set $p =\max(p_1,\ell_{\psi})$. For $i \leq p$, for any $N\in
\bK$, for any $\abar \in N$ define by induction $A^N_0 =A^N_0(\abar)=
\abar$ and $ A^N_{i+1} = A^N_{i+1}(\abar) = \cl^p_N(A_i)$.  Now
applying
Lemma~\ref{compress}, choose $\ell_{\phi}$ so that for every $\abar$
of length $r$, and every semigeneric $N$, $A^N_{p}(\abar) \subseteq
\cl^{\ell_{\phi}}_N(\abar)$.

We want to show that for any semigenerics $M$ and $M'$, for any
$\abar\in M^r$, $\abar'\in M'^r$, and $b \in M$ if
$\cl^{\ell_{\phi}}_M(\abar)
\iso \cl^{\ell_{\phi}}_{M'}(\abar')$ then there is a $b' \in M'$ with
$\cl^{\ell_{\psi}}_M(\abar,b)\iso \cl^{\ell_{\psi}}_{M'}(\abar',b')$
by an isomorphism taking $\abar$ to $\abar'$.
Let
$H_0$ be the substructure of $M$ with universe $(\abar,b)$ and $H_1 =
\cl^{\ell_{\psi}}_M(H_0)$.

Fix $g$ which maps $\abar$ to $\abar'$ and
$\cl^{\ell_{\phi}}_M(\abar)$ isomorphically onto
$\cl^{\ell_{\phi}}_M(\abar')$.
By the choice of $\ell_{\phi}$, for each $i \leq p$, $g$ maps
$A^M_i(\abar)$ isomorphically onto $A^{M'}_i(\abar')$. (Use
Lemma~\ref{1.9a} and induct.)
To avoid superscripts,
for each $i$, let $A'_i$ denote the image of $A_i= A^M_i$ under $g$.
Notice that for some
$j \leq p$, $$(A^M_{j+1}-A^M_j)\inter (H_1-H_0) = \emptyset.$$  Since
$p >
|H_1 - A^M_j|$ this implies
$A^M_j \subm A^M_j H_1$.

Since $M'$ is semigeneric, $M' \sat \phi^p_{A_j,H_1,A_{j+1}}$.
Thus,
there is an isomorphism $\hat g $ extending
$g$ and mapping $H_1$ into $M$ with $
\cl^p_{M'}(A'_j{\hat g}H_1 ) =
\cl^p_{M'}(A'_j)\union{\hat g}H_1$ and so that
$M' |(\cl^p_{M'}(A'_j){\hat g}
H_1)$ is a free join of $\cl^p_{M'}(A'_j)$ and
${\hat g}H_1$
over ${A'_j}$.
Let $H'_1 = {\hat g} H_1$
and $b' = {\hat g}(b)$.  We need to show
$\cl^{\ell_{\psi}}_M(\abar,b)
\iso \cl^{\ell_{\psi}}_{M'}(\abar',b')$.

By the choice of
$\hat g$ and $H'_1$,
$A'_{j}H'_1 \cong A_{j}H_1$ which
contains $\cl^{\ell_{\psi}}_{M}(\abar,b)$,
so it
suffices
(by
Lemma~\ref{1.9a})
to show
$A'_{j}H'_1 $
contains $\cl^{\ell_{\psi}}_{M'}(\abar',b')$.
Note $\cl^{\ell_{\psi}}_{M'}(\abar',b') \subseteq
\cl^p_{M'}(\abar',b') \subseteq
\cl^p_{M'}(A'_j{\hat g}H_1 ) =
A'_{j+1} H'_1$.
By
Lemma~\ref{1.9a},
$\cl^{\ell_{\psi}}_{M'}(\abar',b') =
\cl^{\ell_{\psi}}_{A'_{j+1}H'_1}(\abar',b')$.
Since $A'_{j+1}$ and $H'_1$ are freely joined over $A'_j$,
${\hat g}^{-1}$ is a $1-1$ homomorphism from
$A'_{j+1}H'_1$ onto $A_{j+1}H_1$.  Applying {\bf A5} from
Lemma~\ref{yaxtrue}, we
see $\cl^{\ell_{\psi}}_{A'_{j+1}H'_1}(\abar',b') \subseteq
{A'_{j}H'_1}$ whence $
\cl^{\ell_{\psi}}_{A'_{j+1}H'_1}(\abar',b') =
\cl^{\ell_{\psi}}_{A'_{j}H'_1}(\abar',b')$.


The proof of the following corollary encompasses the derivation of
Theorem~\ref{qethm} from Lemma~\ref{heart}.

\begin{cor}
\label{thmnmc}
Suppose there is a $(\bK_0,\subm)$-semigeneric $L$-structure.  The
theory of the class of $(\bK_0,\subm)$-semigeneric $L$-structures is
nearly model complete.
\end{cor}

\proof.  We have shown that in each semigeneric model the truth of
$\phi(\abar)$ is determined by the isomorphism type of
$\cl^{\ell_\phi}_M(\abar)$ and does not depend on the
particular embedding of $\cl^{\ell_\phi}_M(\abar)$ in $M$. There are
only finitely many
possibilities for this closure and each is determined by a
conjunction
of existential and universal sentences (specifying which $B$ with
$|B|< \ell_{\phi}$ and with $\abar$ enumerating an intrinsic
substructure of $B$ occur).

\begin{cor}
\label{complete}
The theory of  the semigeneric models is complete.
\end{cor}

\Proof.  If $N$ is semigeneric, $N \in \bK$ so, by A6, $N$ does not
contain any substructure $A$, with $\delta(A)< 0$.  Thus,
$\cl_N(\emptyset) = \emptyset$; completeness follows from
Lemma~\ref{heart}.

Recall from \cite{BaldwinShiJapan}:

\jtbdef ${\bK_{0}}$ has the {\em full amalgamation property} if $B
\inter C = A$ and $A \subm B$ implies $D = B \bigotimes_{A} C \in
{\bK_{0}}$ and $C\subm D$.

The following result is proved in \cite{BaldwinShiJapan}

\begin{thm}
\label{talphaexists}
$(\bK_{{\alpha}}, \subm )$ has the full amalgamation property.
There is a generic model $M_{{\alpha}}$ and the theory $T_{{\alpha}}$
of
this generic model is stable.
\end{thm}

Using the full amalgamation property, it is easy to see

\begin{lemm}
\label{genissemigen}
The generic model $M_{{\alpha}}$ for $\bK_{{\alpha}}$ is semigeneric.
\end{lemm}

Combining the above results we have

\begin{thm}  $T_{{\alpha}}$ is nearly model complete.
\end{thm}

The strength of this remark is emphasized by the following
observation.

\begin{thm}  The theory $T_{\alpha}$ is not model complete.
\end{thm}

\proof.
If $T$ is model complete with generic $M$, the type of any finite
subset $X$ is determined
by positive assertions of the substructures that contain $X$.
Fix $A <_i B <_i C \in \bK_0$.  Suppose $A_1, A_2 \subset M$ with
$f:A \iso A_1 $ and $f:A\iso A_2$ and suppose $A_1 <_i B_1 \subm M$,
$A_2 <_i C_2 \subm M$ with 
$f$ and $g$ extending to $\hat f$, $\hat g$ such that ${\hat g}:B  
\iso B_1$ and ${\hat g}: C\iso C_2$.   Then every
existential formula true of $A_1$ is true of $A_2$ but the converse
is obviously false.   To see the nonobvious assertion, let $D$ be  
arbitrary with $A_1 \subseteq D \subset M$ and $D$ not contained in
$B_1$.  Then,  $B_1 \subm B_1D$.  Pulling back to $A,B,C$, by full
amalgamation,  there is 
a $D' \supset D$ such that $BD' \otimes_{B} C \in \bK_0$ and
$C \subm CD'$. Extending ${\hat g}$ from $C$ to $CD'$ provides the  
required witness.

\jtbnumpar{Remark} Let $L$ contain a single binary relation and
restrict to the class of graphs.  Baldwin and Shi noted
\cite{BaldwinShiJapan} that
full amalgamation holds for the class $\bK'_{{\alpha}}$ consisting of
those graphs in $\bK_{{\alpha}}$ which omit squares.  Laskowski
observed
that this argument applies as well to the class $K^n_{{\alpha}}$ of
graphs
which
omit cliques of size $n$.  Thus the theory of the generic model
associated with each of
these classes is stable and nearly model complete.

\section{Adding Random relations}
\label{addran}

In this section we begin with the collection of finite models for a
language $L$ with only the equality symbol (i.e. $n$-element sets for
arbitrary $n$)
and add
additional `random' relations with respect to probability measures
described below.

  We show that
a $0-1$ law holds for the set of first order sentences in the
expanded language
and that the almost sure theory (the sentences with limit probability
$1$)
 is stable.
Adding a single symmetric irreflexive binary relation
gives the family of theories investigated independently by
Shelah-Spencer
and Baldwin-Shi.
Viewing this situation as an expansion of the language of equality
may
seem eccentric but we expect to exploit this viewpoint for more
interesting base languages in the future.  This project is
well-advanced in
\cite{Shelah550}.

\jtbnumpar{Context}
\label{exa}
Let $L$ contain only the equality symbol.  The $L$-structure $M_n$ is
a set with $n$ elements.
$\bK_0$ is the class of all finite sets and $\bK$ the class of all
sets.  On $\bK$, $\subm$ is just $\subseteq$ and $A \subi B$ just if
$A= B$.

\jtbnumpar{Remark}
The properties {\bf A1}-{\bf A6} and the conclusion of Lemma~\ref{h1}
hold for $\bK$ in Context~\ref{exa}. Moreover, $\bK_0$ has the full
amalgamation
property.

\jtbdef
\label{primdef}
We say that $B$
is a {\em primitive} extension of $A$ if $A \subm B$
and for every $B'$ with $B'$ properly contained between $A$ and $B$,
$B'$ is not a strong submodel of $B$.

Now, we show how to define the notion of independent random relations
(with edge probability `$n^{-\alpha}$') for an
arbitrary finite relational language $L^+$.  Then we define the
notions
of
dimension and strong submodel in the  extended language $L^+$ and
show
that the properties {\bf A1}-{\bf A5}  hold for the extended language
and {\bf A6} holds with probability 1.

\jtbnot We write $[X]^m$ for the collection of $m$-element subsets of
a set $X$.  We will write either $C \in [X]^m$ or  (surreptitiously
fixing an enumeration of $C$) $\cbar \in [X]^m$ to indicate a member
of this set.

\jtbnumpar{Adding Random relations}

\label{probdef}

Fix an enumeration $\langle R_i: i<p\rangle$ of the relation symbols
in $L^+-L$ and let $k_i$  denote the arity of $R_i$.  Let $L_i$
contain only $R_i$.
Let $t $ denote the largest arity of the $R_i$.  Fix also  a sequence
of numbers $\alpha_i$ with $0< \alpha_i \leq 1$  and $\gamma_i$ with
$0 \leq \gamma_i \leq 1$  for $i<p$.  (We will require later that the
$\alpha_i$ and $1$ be linearly independent over the rationals.)

We will define for each isomorphism type of an $L^+$ structure of
size $n$, the probability of a random structure of size $n$, having
that isomorphism type.

We  assume that each new relation in
the expanded structure is symmetric and irreflexive in the sense of
Paragraph~\ref{fixnotation}.
Note that this
formalism does not describe what one should mean
by a
random directed graph.

Let $N$ be an $L^+$ structure of cardinality $n$.  Let $C$,
enumerated as $\cbar$, be a subset of $N$ with size $k_i$.  Let
$$q_{i,n}(C|L_i) = \left\{
\begin{array}{ll}
\gamma_i n^{-\alpha_i}&
\mbox{ if }
N\sat R_i(\cbar)\\
1-\gamma_i n^{-\alpha_i}&\mbox{ if  }N\sat \neg R_i(\cbar)
\end{array}
\right.    $$

and,
$$P_n(N) = \prod_{i < p}\prod\{q_{i,n}(C): C \in [N]^{k_i}\}.$$

\jtbnumpar{Remark}
\label{unpack}
Let $w_j(N) = |\{C \in [N]^{k_j}: N \sat R_j(\cbar)\}|$ and\\
$\wbar_j(N) = |\{C \in [N]^{k_j}: N \sat \neg R_j(\cbar)\}|$.  Then

$$
P_n(N)= \prod_{j < p}(\gamma_j n^{-\alpha_j})^{w_j(N)}
(1-\gamma_jn^{-\alpha_j} )^{\wbar_j(N)}.
$$

If $L^+$ has a single binary edge relation and the probability of a
two
element structure is $n^{-\alpha}$ when the points are related
and $1-n^{-\alpha}$ if not, we
return to the situation of \cite{ShelahSpencer}.

Recall from Lemma~\ref{axsemigeneric} the sentences axiomatizing
the class of semigeneric models.
We want to show that the almost sure theory exists and is exactly the
theory of the
semigenerics.  To this end, we will show $$\lim_{n \rightarrow
\infty}P_n(\phi^m_{A,B,C})=1$$ for each $m, A, B, C$.

\jtbnot  Henceforth, $A, B, \ldots M,N \ldots$ range over
$L$-structures.
$A^+, B^+$ etc. denote an expansion of $A$, respectively $B$ to
$L^+$.
We refer to
the universe of $A^+$ or $A$ by either of these terms rather than the
more accurate $|A^+|$ or $|A|$ and reserve $|\,\, |$ for cardinality.
Thus $A^+|L = A$ and we use these notations interchangeably.

We now translate our probability asssignment into a class
$(\bK,\delta)$ as in
Example~\ref{kalpha}.

\jtbnot
\label{newdim} Let $\bK^*_0$ be the collection of all finite $L^+$
structures.
\begin{enumerate}
\item For $A^+ \in \bK_0^*$, define $\delta(A^+)
=\delta_{\alpha}(A^+)$  as in
Example~\ref{kalpha}, using only the relation symbols in $L^+-L$ and
using the
parameters $\alpha_i$ from Paragraph~\ref{probdef}.
\item $e(B^+/A^+)$ denotes $e(A^+B^+) - e(A^+)$.
\item  $\bK^+_0$ denotes the collection of $A^+ \in \bK^*_0$ such
that for each
$A' \subseteq A$, $\delta(A') \geq  0$.
\item $\gamma(A^+,B^+) = \prod_{i< p}  \{\gamma_i:|C| = k_i, C
\subseteq B,C \not \subseteq A \& B \sat R_i(\cbar)\}$.
\end{enumerate}

\jtbnumpar{Remark}
\label{fixinga6}  The link between the function $\delta$ and the
probabilistic situation
is provided in Remark~\ref{justdelta} where we show that the
expectation of the existence
of a copy over $A^+$ of a structure $B^+$ is determined by
$\delta(B^+/A^+)$.
Lemma~\ref{justa6}  implies that
replacing  $\bK_0^*$ by  $\bK_0^+$ does not change the almost sure
theory.  But it does make the model
theory conform with the framework of Section 1.  The major
calculations of this section
were carried out in $K^*_0$ in \cite{ShelahSpencer}.

\jtbnumpar{Assumption}[Irrationality Hypothesis] The coefficients
$\alpha_i$ and $1$ are  linearly independent over the rationals.
This
generalizes for an arbitrary finite language the assumption
in the case of random graphs with $p = n^{-\alpha}$ that
$\alpha$ is irrational.  This hypothesis easily implies that for any
$L$-structures
$A \subm B$ and expansions
$A^+ \subseteq B^+$, $\delta(B^+/A^+)\neq 0$.  This is a key property
(see Remark~\ref{reasonforirration}).

We have the notions of $\subm$ and $\subi$ as a relation on members
of
$\bK_0$. The definition of $\delta$ on $\bK^+_0$ induces  
corresponding 

relations on
$\bK^+_0$.  Since
we work directly with $\delta$ it is worthwhile recording the
translation.

\jtbdef
For $A^+\subseteq B^+ \in \bK^+_0$,
\begin{enumerate}
\item $A^+ \subm B^+$ if
$\delta(B^+_1/A^+) > 0$ for every $B^+_1$
with
$A^+ \subset B^+_1 \subseteq B^+$.
\item
$A^+ \subi B^+$ if
$\delta(B^+/B^+_1) < 0$ for every
$B_1^+$ with $A^+\subseteq B^+_1 \subset B^+$.
\item
$B^+$ is a {\em primitive extension} of $A^+$, if $\delta(B^+/A^+ )>
0$ and
$\delta(B^+/A_1^+) \leq 0$ for each $A_1^+$ with $A^+ \subset A_1^+
\subseteq
B^+$.
\end{enumerate}

Note that:

\jtbnumpar{Remark}  $\bK^+_0$ satisfies axioms {\bf A0}-{\bf A6}.
\label{k+good}
We will be using the following monotonicity properties which follow
formally as in Section 1.

\begin{lemm}
\begin{enumerate}
\item
If $A^+ \subi B^+$ and $A^+ \subseteq B_1^+ \subseteq B^+$
then $B^+_1 \subi B^+$.
\item
If $A^+ \subm B^+$ and $A^+ \subseteq B_1^+ \subseteq B^+$
then $A^+ \subm B^+_1$.
\end{enumerate}
\end{lemm}

\jtbnumpar{Remark}
\label{needmn}
The exact phrasing of the following notions is extremely delicate.
We
consider a fixed pair of finite $L^+$-structures, $A^+ \subm B^+$.
The
$L$-structures $\langle M_n:n < \omega \rangle$ naturally form a
chain
so an embedding $f$ of $A$ into $M_n$ can naturally be regarded as a
map
of $A$ into $M_m$ for $m > n$.  We are concerned with the properties
of
{\em extensions} of $f$.  Thus,  the immediately following definition
of an
$L^+$-homomorphism extending $f$ is agnostic concerning the
preservation of
relations on $A$.

\jtbdef
\label{fixranvar}
Let $A^+ \subseteq B^+$.
Let $f$ be a $1-1$ map from  $A$ into $M_n$, and let $G^+$ be
an $L^+$ structure  expanding $M_n$.
Let $T$ denote the range of $f$.
\begin{enumerate}
\item
We say
an injective map $g:B^+ \mapsto G^+$ which extends $f$ is an
{\em $L^+$-homomorphism} relative to $A$is an
if for  any $L^+$-relation $R$,  and any $\bbar\in
B$ but
not in $A$, $B^+ \sat R(\bbar)$
implies $G^+\sat R(g\bbar)$.
 \item
For any $G^+$ expanding $M_n$, and $W \subseteq M_n$ with $|W| \leq
n$  we say $N(f,A^+,B^+,W)  = k$ in $G^+$, if

$$k =  |\{g:B^+ \mapsto W \subseteq M_n  \text{ is an
$L^+$-homomorphism relative to $A$ and } g \supseteq f
\}|.$$
If $W = M_n$ we omit it.
\item
We say $G^+$ is in the event $Y_{f}$,
which depends on a constant $c_1$,  if
$$
n^{\delta(B^+/A^+)} (\log n)^{-(v(B/A)+1)} < N(f,A^+,B^+) < c_1
n^{\delta(B^+/A^+)}.    $$
\item
Let $U$ denote the range of $f$.
For each $S \subseteq M_n$ with $S
\inter U  = \emptyset$ and $|S|
= |B-A|$, fix
(if possible)
an $L$-isomorphism $g_S$ between $B$ and
$US$ which extends $f$.  (Since $L$-isomorphism just means $1-1$ map,
such $S$ and $g_S$ exist whenever $n \geq |B|$).
\item
  For each such $S$ with fixed
$L$-isomorphism $g_S$ of $B$
into $M_n$,
let $X_{f,S}$ be a random variable such that $X_{f,S}(G^+)$
is
$$\left\{
\begin{array}{ll}
1&
\mbox{if  $g_S$ maps $B$ onto $US$ and is an
$L^+$-homomorphism relative to $A$ into $G^+$}\\
0&\mbox{ otherwise.  }
\end{array}
\right.    $$

For $W \subseteq M_n$ with $|W| \leq n$,
let
$$X_{f,W} = \sum \{X_{f,S}: g_S:B \mapsto M_n \text{ and }
g_S\supseteq f \text { and } S \subseteq W\}.$$
If $W = M_n$, we write $X_f$ for $X_{f,W}$.
\end{enumerate}

\jtbnot
\begin{enumerate}
\item
 For any property $P$ of  structures,  in particular a first
order property, the assertion, `for almost all sufficiently large
$M$,
$M \sat P$' (abbreviated a.a.)
means `for every $\epsilon >0$ there is an $N$ such that if
$n > N$, $P_n(\{G^+|L =M_n: G^+ \sat P\}) >1-\epsilon$'.
\item By an {\em indicator random variable} we mean one which takes
values $0$ or $1$ and thus indicates a set.
\item We write $f \iso g$ if $f = O(g)$ and
$g = O(f)$.
\end{enumerate}

The next  lemma expresses the key observation linking the
probability with the dimension function $\delta$.

\begin{lemm}
\label{expect}
For all sufficiently large $n$
and all $f:A \mapsto M_n$,
and any $W \subseteq M_n$ with $|W| \leq n$, the expectation
$$
E(X_{f,W})
\iso |W|^{v(B/A)}
n^{-e(B^+/A^+)}.$$
\end{lemm}

\Proof.  The probability of an $L$-embedding of $B^+|L$ into $M^+$
actually being an $L^+$ homomorphism is
$\gamma(A^+,B^+)n^{-e(B^+/A^+)}$.
The number of such embeddings has order of
magnitude $|W|^
{v(B/A)}$.
Since expectation is additive this yields
$$
E(X_{f,W})
\iso
|W|^{v(B/A)}
\gamma(A^+,B^+)
n^{-e(B^+/A^+)}$$

The constant is absorbed by the approximation $\iso$.  In particular,
we have:

\jtbnumpar{Remark}
\label{justdelta}

If $W = M_n$ this simplifies to
$$
E(X_f)
\iso
n^{\delta(B^+/A^+)}.$$

In Theorem~\ref{semigen} we guarantee that our extensions are
$L^+$-isomorphisms (no
new relations) rather than just $L^+$-homomorphisms.
Now we justify the restriction from $\bK^*_0$ to $\bK^+_0$.

\begin{lemm}
\label{justa6}
If $\delta(B^+) <0$ then a.a.~there is no embedding of
$B^+$ into $G^+$.
\end{lemm}

\proof.  The expected number of copies of $B^+$ is $n^{\delta(B^+)}$.
If $\delta(B^+) < 0$, this tends to $0$.

\begin{thm}
\label{bounding}

Fix $A^+ \subseteq B^+$ with
$A^+ \subm B^+$.    

Let $V$ be the event: for all $
f:A\mapsto M_n$, the event $Y_f$ holds. 

 Then, for some choice of $c_1$ (recall $Y_f$ depends
on $c_1$.),

$$\lim_{n\rightarrow \infty}P_n(V)
=1.$$
\end{thm}



\proof.  By a straightforward induction, we can reduce to the case
that
$B^+$ is a primitive extension of $A^+$.
The proof of this case proceeds through several definitions and
lemmas.
Considering the definition of $Y_f$, one can see that we need to
establish both lower and upper bounds.  The lower bound argument
proceeds as follows.
Roughly speaking, for $f:A \mapsto M_n$ and $W \subseteq M_n$, we say
$(f,W)$ is bad if there is no extension of $f$ to an
$L^+$-homomorphism
(in the sense of Definition~\ref{fixranvar})
of $B^+$ into $W$.  In Lemma~\ref{firstpass} we show that if
$W$ meets a cardinality requirement specified in
Definition~\ref{goodW}
then the probability that $(f,W)$ is bad is less than 1/2.  By
strengthening the requirements on $W$ as in Definition~\ref{betterW}
we improve the upper bound on the probability that
$(f,W)$ is bad in Lemma~\ref{secondpass}.  Finally, taking into
account
the number of possible $W$'s, we complete the proof of the lower
bound
in Paragraph~\ref{lower}.
After several preliminary definitions and lemmas we complete the
proof
of the upper bound in Paragraph~\ref{upper}.

\jtbdef
\label{goodW}
We have fixed $A^+ \subseteq B^+$ with
$B^+$ a primitive extension of $A^+$.
  For $G^+$ an $L^+$-expansion of $M_n$
and $W \subseteq M_n$
and $f$ an $L$-isomorphism of $A$ into $M_n$,
$(f,W)$ is {\it bad} in $G^+$ if there is no $g$ defined on $B-A$
into $W$
such that $f \cup g$ defines an $L^+$-homomorphism from $B^+$ into
$G^+$.

\begin{lemm}
\label{firstpass}
There is a constant $s$ such that for all sufficiently large $n$ and
any
$L$-isomorphism
$f:A \mapsto M_n$, if $|W| $ is the integer $m =m_s =
[sn^{e(B^+/A^+)/v(B/A)}]$
then $$P_n((f,W) \text{ is bad })< 1/2.$$ \end{lemm}

\proof.  Without serious loss of precision, $W \cap \rng f =
\emptyset$.
We use the notation from Definition~\ref{fixranvar}.

Now $$P_n((f,W) \text{ bad }) = P_n (X_{f,W} = 0)$$
so we want to show that for all sufficiently
large $n$,
$$P_n(X_{f,W}=0)< 1/2.$$

By Chebyshev's inequality,
$$P_n(X=0) \leq {\Var (X) \over E(X)^2}.$$

By Lemma~\ref{expect},
$$E(X)=
E(X_{f,W})
\iso
|W|^{v(B/A)} n^{-e(B^+/A^+)}.$$
Using the fact that $|W| = m_s$ this shows
$E(X)$ is a polynomial of degree $v$ in $s$, as the powers of $n$
cancel. We will obtain the required result by showing $\Var(X)$
is a polynomial of degree $2v-1$ in $s$
which implies that for sufficiently large $s$,
${\Var (X) \over E(X)^2}
< 1/2$.

Now,
$$\Var(X) = \sum_S \Var(X_S) + \sum_{S \neq T} \Cov(X_S,X_T)$$
where $S,T$ range over subsets of $M_n$ disjoint from the image of
$f$.

An easy calculation shows that for
any set of  indicator random variables $$E(\sum_S X_S) \geq
\sum_S\Var(X_S);$$ so we have

$$\Var(X) \leq E(X) + \sum_{j=0}^{v}
\sum_{
|(S \inter T) -A| =j}
\Cov(X_S,X_T).$$

If
$|(S \inter T) -A| =0$
then
$S \inter T = \emptyset$ and $\Cov(X_S,X_T)$ is zero.
 Always, $\Cov(X_S,X_T)\leq E(X_S X_T)$ which,
since these are indicator random variables, is just
$P_n(X_SX_T)$.  Recall the definition of the probability measure
from Definition~\ref{probdef}; $t$ is largest arity in the language.
$$P_n(X_SX_T)=
\prod_{m \leq t}
\prod_{A \in [ST]^{m}}q_m((M|L^{\leq m})|A)$$

$$
\leq
\prod_{m \leq t}({{
\prod_{A \in [S]^{m}}q_m((M|L^{\leq m})|A)
\prod_{A \in [T]^{m}}q_m((M|L^{\leq m})|A)}
\over{
\prod_{A \in [S\inter T]^{m}}q_m((M|L^{\leq m})|A)}}.$$

Let
$B'$ be $g_s^{-1}(S \cap T)$.  So $|B'|=j$.
Abbreviating the notations from \ref{newdim}, let $c_S = c_T =
\gamma(A^+,B^+)$ and $c_{B'} = \gamma(A^+,B')$.  Similarly, let
$e=e_S = e_T =
e(B^+/A^+)$ and $u' = e(B'/A^+)$.
With this notation we can rewrite the last inequality as
$$P_n(X_SX_T)\leq {{c_S c_T }\over {c_{B'}}}n^{-(2e-u')}.$$
(The key to the inequality is that $2e-u'$ may undercount the
number of relations on $ST$ but this undercount can only overestimate
the probability).
$|(S \inter T) -A| =j$, $|ST-A| = 2v -j$ so
$\chi_W(ST/A) \iso m^{2v-j}$.

If ${u'\over j} >    {e\over v}$ then ${e\over v} > {{e-u'}\over
{v-j}}$
which contradicts the fact that
     $B^+$ is a primitive
 extension of $A^+$.                  So  $u' \leq {je\over v}$.
Thus,
$$
\sum_
{|(S \inter T) -A| =j}
\Cov(X_S,X_T)
\leq m^{2v-j}n^{u'-2e}\leq (sn^{e/v})^{2v-j}
n^{je/v-2e} =
{(s^v)}^{2-j/v}.$$
(We can drop the constants in the last computation as ${c_Sc_T \over
c_B'}< 1$.) So,
$$\Var(X) \leq E(X) + \sum_{j=1}^{v}
(s^v)^{2-j/v}\leq  E(X) + vs^{2v-1}.$$
Since $E(X)$ has degree $v$ in $s$, this implies
$\Var(X) \leq  E^2(X)/2$ for sufficiently large $s$ and so
 $$P_n((f,W) \text{ is bad })< 1/2.$$
\medskip

Now we want to modify the choice of $W$ to get a better upper bound
on
the probability that $(f,W)$ is bad.

\jtbdef
\label{betterW}
Choose $s$ by Lemma~\ref{firstpass}.  As before,
let $m_s =[s n^{e(B^+/A^+)/v(B/A)}]$.
We say that $W \subseteq M_n$ is {\em $k$-appropriate} if
$|W|=
  [1 +k{m_s} \ln n] $.

\begin{lemm}
\label{secondpass}
For all sufficiently large $n$ and any
$L$-isomorphism
$f:A \mapsto M_n$,   

 for sufficently large $k$, if $W \subseteq n$ is $k$-appropriate,
$$P_n((f,W) \text{ is bad }
)< {{n^{-|A|-1}\over 2}}.$$
\end{lemm}

\Proof.  Again, assume
without loss of generality that $W \cap A = \emptyset$.
Suppose
$W$ contains
$k\ln n$ disjoint subsets $W_i$ each with cardinality ${m_s}$.
For $(f,W)$ to be bad, each of the $k\ln n$
independent events that $(f,W_i)$ is bad must occur and by
Lemma~\ref{firstpass} $P_n((f,W) \text{ is bad })<1/2$.  Thus,

$$P_n((f,W) \text{ is bad })<
2^{-k\ln n}.$$
But for all sufficiently large $n$ and $k$,
$$2^{-k\ln n}
< {{n^{-|A|-1}\over 2
}}$$
so we have the result.
\medskip

We have shown that for each $f$, a.a. there is a $W$ such that
$(f,W)$
is not bad.  The next paragraph strengthens this assertion.
\jtbnumpar{Proof of Lower Bound in Theorem~\ref{bounding}}
\label{lower}

Fix $k$ satisfying the conclusion of Lemma~\ref{secondpass}.
For an $L$-isomorphism  $f$ of $A$
into $M_n$,
let the random variable $Z_f(G^+)$
be the number of  $k$-appropriate $W \subseteq M_n$( i.e.$|W|= w =
  [1 +k{m_s} \ln n]
$ with $m$ from Definition~\ref{betterW})  such that $(f,W)$ is bad.
Let $\gamma$ denote the number of possible $k$-appropriate
$W$.
(The value of $\gamma$ is not used
in the first stage of the argument.)
Then,
$E(Z_f)< {\gamma}n^{-|A|-1}$. So, by Markov's inequality,
$$P_n(Z_f \geq {\textstyle{\gamma}/2)} \leq 2 E(Z_f)/{\textstyle
{\gamma}} <
2n^{-|A|-1}.$$
 But then, since there are only $n^{|A|}$ choices for $A$,
a.a.~for {each} $f$ at most half of the $W$ are bad for $f$.

Let $v$ denote $|B-A| = v(B/A)$.
Each extension $g$ of $f$ to $B$ is contained in at most
$
 {{n-v} \choose {w-v}}$
$k$-appropriate $W$, since there are approximately
$ {{n-v} \choose {w - v}}$ choices for the  elements
which comprise $W - \im g$.
So at most
{$|N(f,A^+,B^+)|
{n-v \choose {w -v}}$}
{$k$}-appropriate $W$
contain an extension of $f$  but at least $\gamma/2$ do.
Now note that $ \gamma = {n \choose {w}}$

Thus, a.a.~for all $f$,
$$|N(f,A^+,B^+)|
{n-v \choose
w-v} \geq
{{1 \over 2} \displaystyle{n \choose {w}}}.$$
Noting that   ${{{n-v} \choose {w- v}}}
$ is approximately ${{n \choose w}({w\over n})^{v}}$,
 we have
$$|N(f,A^+,B^+)| \geq
n^{-e(B^+/A^+)}
(n/ w)^{v}/2.$$
Recalling that $w = [1 +k{m_s} \ln n]$, this implies for every $f$,  
a.a.
$$|N(f,A^+,B^+)|
> n^{\delta(B^+/A^+)}(2ks\ln n)^{-v}$$
which establishes the lower bound $n^{\delta(B^+/A^+)}(\ln n)^{-c}$  
by taking
$c = v+1$.

\jtbnumpar{Remark}  The statement and proofs of
of the probability analysis are based on the argument
in \cite{ShelahSpencer}.
The first author
acknowledges discussions with Albert, Cherlin, Lachlan, and Laskowski
on the details of
the current argument, and
supplemental remarks to the original paper by Spencer.

 \jtbnumpar{Remark}
\label{reasonforirration}
  The irrationality hypothesis
is necessary to make fruitful application
of this result.  If
there exist $A\subm B$ with $\delta(B/A) = 0$
then the lower bound we have
established is less than one rather than tending to infinity
as $n$ does.  This destroys the argument of Theorem~\ref{semigen}.

\jtbnumpar{Remark}
\label{fewins}  From Lemma~\ref{h1} we have:
  Let $A^+\subi C^+$.  There exists a $K$
such that a.a.~for every embedding $f$ of
$A^+$ into an expansion $G^+$ of $M_n$, there are fewer than $K$
{$L^+$}-homomorphisms extending $f$ from $C^+$ into $G^+$.

\medskip
\jtbnumpar{Proof of upper bound in Theorem~\ref{bounding}}
\label{upper}
Since $B^+$ is a primitive extension of $A^+$, $(Ab)^+ \subi B^+$ for
any $b\in B-A$.  Thus, by
Remark~\ref{fewins} there are fewer than $K$ extensions $g_i$ with
any fixed image of $b$.
The range of each extension $g_i$
can intersect at most $K |B|^2$ other extensions
so if $N(f,A^+,B^+) =s'$, there is a set of $s = s'/(K |B|^2+1)$
disjoint
extensions. Let $p = |A|$, $v = v(B/A)= |B-A|$ and $e= e(B^+/A^+)$.
For an appropriate constant $c<1$,
there are less than
$cn^{p}(n^{vs}/s!)$
pairs of a function $f$ taking $A$ into $M_n$
and a set of $s$ extensions
(disjoint over $\rng f$)
$\langle g_1, \ldots g_s\rangle$.  The probability that
each of the
$g_i$ is an $L^+$-homomorphism is $n^{-e}$ so the
probability of such a pair of a function and $s$
homomorphisms is at most
$n^p(n^{vs}n^{-es})/s!$
By Stirling's formula, this is less than
${\displaystyle {n^p(2.72 n^{v})^sn^{-es}}\over
{(2\pi)^{1/2} s^{s+1/2} (2.72)^{1/(12 s+1)}}}$
which is much less than
$1$ if $s \geq 3n^{\delta(B^+/A^+)}$.  (
Observe that $n^p({2.72\over 3})^s$ tends to $0$ as $n$ tends to
infinity.)
 Thus, a.a.~$s\leq 3n^{\delta(B^+/A^+)}$. 

a.a.~for each f,
$$N(f,A^+,B^+) = s' \leq 3(K|B|^2+1)n^{\delta(B^+/A^+)}$$ proving
Theorem~\ref{bounding}.

\medskip
We now want to show that each of the axioms for semigenericity has
limit
probability $1$.  Roughly, the program is to show that for $A^+ \subm
B^+$ and $f:A\mapsto M_n$, the number of extensions of $f$ to
$1-1$-homomorphisms of
$B^+$  is much greater than the number of such extensions which fail
to
witness
the
definition of semigenericity.  Since there are a bounded number of
types of failure, it suffices to check each type separately as we
do
in the following argument.

In general, embeddings $f:A^+ \mapsto M_n$ and ${\hat f}:B^+ \mapsto
M_n$ fail to witness semigenericity of $G^+$ if
\begin{enumerate}
\item $\cl^m_{G^+}({\hat f}B) \neq {\hat f} B \cup \cl^m_{G^+}(fA)$
or
\item
$\cl^m_{G^+}(fA)$ and
${\hat f} B$ are not freely joined over
${fA}$ in $G^+$.
\end{enumerate}

In considering $\phi^m_{A^+,B^+,C^+}$, we are fixing on $C^+$ as a
specific candidate
for the isomorphism type of $\cl^m_{G^+}(f A)$.

\begin{thm}
\label{semigen}
  If $A^+\subm B^+$ and $A^+ \subi C^+$ with $|C^+| < m$ then
$$\lim_{n \rightarrow \infty}P_n(\phi^m_{A^+,B^+,C^+})=1.$$
\end{thm}

\proof.
For any $f$ mapping $C$ into $M_n$,
and a $1-1$ homorphism $f'$ extending $f$ to $E$,
$(G^+,f')$ fails as a witness
for $C$ and $f$ if
\begin{enumerate}
\item $f'$ is not an $L^+$-isomorphism or
\item $\cl^m_{G^+}({ f'}B) \neq { f'} B \cup fC$ or
\item
$fC$ and $f'B$ are not freely joined over
${fA}$ in $G^+$.
\end{enumerate}

Note that $C^+\subm C^+ \otimes_{A^+} B^+$ and, letting $E^+$
denote $C^+ \otimes_{A^+} B^+$,
$\delta(B^+/A^+ )=  \delta(E^+/C^+)$.
By Theorem~\ref{bounding}, more specifically Paragraph~\ref{lower},
a.a.~for each $f$,  $$|N(f,C^+,E^+)| >n^{\delta(B^+/A^+)}
(\log n)^{-c}$$
 where $c = v(B^+/A^+) +1$.
For conditions i) and  iii) consider
any $F^+$ which is an expansion of $E^+$ by adding additional
relations.
Then $\delta(F^+/A^+)= v(E/A) -
e(F^+/A^+)$
and
$e(F^+/A^+)
>e(E^+/A^+)$.
By Theorem~\ref{bounding}
$$N(f,C^+,F^+) < c_1 n^{\delta(F^+/A^+)}<n^{\delta(E^+/A^+)} (\log
n)^{-c}.$$

For condition ii) for any $D^+ \in \Dscr^m_{B^+,C^+}$,
$$\delta(D^+B^+/C^+)<
\delta(B^+C^+/C^+)
=\delta(B^+/A^+).$$
If $C^+$ is not
strong in $D^+B^+$ then by Lemma~\ref{fewins},
$N(f',C^+,D^+B^+) < K$.
If $C^+ \subm D^+B^+$ then by Theorem~\ref{bounding}, more
specifically
Paragraph~\ref{upper},
a.a. $$N(f,C^+,D^+B^+) < c_1
n^{\delta(D^+B^+/C^+)}<
n^{
\delta(B^+C^+/C^+)} (\log n)^{-c}.$$

Now the number of isomorphism types of extensions $C$ that
have failures $f'$ is bounded in
terms of the cardinality
of $A^+$, $B^+$, and $m$; it does not depend on $n$.  If this number
is
$L$,  the
total number
of failures of any sort is less
than $Ln^{\delta(E^+/A^+)} (\log n)^{-c}$.
Thus, the probability that for each $f$, one of
the extensions of $f$ witnesses $\phi^m_{A^+,B^+,C^+}$ tends to one
as
required.

\begin{lemm}
\label{emptycl}
For every $m$,
a.a.~$\cl^m_{L^+}(\emptyset)=\emptyset$.
\end{lemm}

\proof.  $A^+ \subseteq \cl^m_{L^+}(\emptyset)$ just if
$\delta(A^+)<0$. But in passing from $\bK^*_0$ to $\bK^+_0$ (cf.
Lemma~\ref{justa6},
we have forbidden such $A$.

We collect our results in the following theorem which requires the
definition of two theories.

\jtbnot
\label{conclude}
Let $L$ contain only the equality symbol and let
$L^+$ be an arbitrary finite relational language containing $L$.
Suppose
probabilities are defined on finite $L^+$ structures as in
Definition~\ref{probdef} with the $\alpha_i$  and $1$ linearly
independent
over the rationals.
By $T^{\alpha}$, the almost sure theory of random $L^+$-structures
we mean the collection of $L^+$-sentences which have limit
probability
1.    Recall that  $T_{\alpha}$ is the theory of the generic
structures
for $\bK_{\alpha}$ (Definition~\ref{kalpha})
whose existence is guaranteed by
Theorem~\ref{talphaexists}.

A theory $T$ has
the finite model property if every theorem of $T$ has a finite model.

\begin{thm}  Under the hypotheses in Notation~\ref{conclude},
$T^{\alpha}$, the almost sure theory of random
$L^+$-structures is the same as the theory $T_{\alpha}$
of the $\bK_{\alpha}$-generic model.
This theory is
complete, stable, and nearly model complete.  Moreover, it has the
finite model property and has only infinite models
so is not finitely
axiomatizable.
\end{thm}

\Proof.  By Theorem~\ref{semigen} and the choice of $\bK^+_0$, every
model of $T^{\alpha}$ is
$(\bK_{\alpha}, \leq_{\alpha})$-semigeneric.  By
Corollary~\ref{thmnmc},
$T^{{\alpha}}$ is nearly model complete.  By Corollary~\ref{complete}
and Lemma~\ref{emptycl}, $T^{{\alpha}}$ is complete.
Since the generic model for
$\bK_{\alpha}$ is semigeneric, $T^{{\alpha}} = T_{\alpha}$.
\cite{BaldwinShiJapan} shows that $T_{{\alpha}}$ is stable.

Since each theorem of  $T^{{\alpha}}$ has limit probability $1$, for
arbitrarily large $n$, there is nonzero probability that there is a
model of size $n$.  Thus,  $T^{{\alpha}}$ has the finite model
property.

\jtbnumpar{Remark}
The major novelty of this result
is the identification of the two theories, thereby obtaining the
stability of $T^{{\alpha}}$ and the non-finite axiomatizability
of $T_{{\alpha}}$.  The notion of near model completeness specifies
the precise degree of quantifier elimination in $T^{\alpha}$.
In addition, we have extended the $0-1$
law from a language with a single binary relation to
an arbitrary
finite relational language.


\begin{thebibliography}{1}

\bibitem{BaldwinShiJapan}
J.T. Baldwin and Niandong Shi.
\newblock Stable generic structures.
\newblock {\em Annals of Pure and Applied Logic}, 199x.
\newblock to appear.

\bibitem{Bgroup}
A.~Baudisch.
\newblock A new $\aleph_1$-categorical pure group.
\newblock 1992.

\bibitem{Hrustableplane}
E.~Hrushovski.
\newblock A stable $\aleph_0$-categorical pseudoplane.
\newblock preprint, 1988.

\bibitem{KuekerLas}
D.W. Kueker and C.~Laskowski.
\newblock On generic structures.
\newblock {\em Notre Dame Journal of Formal Logic}, 33:175--183,  
1992.

\bibitem{Lynch1}
J.~Lynch.
\newblock Probabilities of sentences about very sparse random graphs.
\newblock {\em Random {S}tructures and {A}lgorithms}, 3:33--53, 1992.

\bibitem{Shelah550}
S.~Shelah.
\newblock 0-1 laws.
\newblock preprint 550, 199?

\bibitem{Shelah467}
S.~Shelah.
\newblock Zero-one laws with probability varying with decaying  
distance.
\newblock Shelah 467, 199x.

\bibitem{ShelahSpencer}
S.~Shelah and J.~Spencer.
\newblock Zero-one laws for sparse random graphs.
\newblock {\em Journal of A.M.S.}, 1:97--115, 1988.

\bibitem{Wagnerdim}
F.~Wagner.
\newblock Relational structures and dimensions.
\newblock In {\em Automorphisms of {f}irst {o}rder {s}tructures},  
pages
  153--181. Clarendon Press, Oxford, 1994.

\end{thebibliography}
\end{document}